\documentclass{amsart}

\usepackage{amssymb,amsmath,amsthm,amscd}%--pas necessaire avec amsart

\newcommand{\N}{\mathbb{N}}%------------------------Integers >=0
\newcommand{\C}{\mathbb{C}}%------------------------Complex numbers
\newtheorem{theo}{Theorem}[section]%
\newtheorem{cor}[theo]{Corollaire}%
\newtheorem{rem}[theo]{Remarque}%

\begin{document}

\title {Courants du type r\'esiduel attach\'es \`a une intersection compl\`ete
}

\author[E.Mazzilli]%
{Emmanuel Mazzilli}%

%\date{}%
%\subjclass[2000]{Primary 53D40, 53D12; Secondary 37D15.}

\address{E.M.: UFR de Math\'ematiques\newline%
 \indent Universit\'e de Lille 1\newline%
 \indent 59655 Villeneuve d'Ascq\newline%
 \indent France} %
\email{mazzilli@math.univ-lille1.fr}

 \keywords{Residues currents, local theory of residues.}%
 \subjclass{32A27, 32A55, 32C30.}%

\begin{abstract}\hskip5pt
We construct in complete intersection's case, elementary currents
which describe the local ideal, and give a decomposition in it for
holomorphic function.
\end{abstract}

\maketitle

\section{Introduction et pr\'eliminaires.}

Si $f$ est une fonction holomorphe sur une vari\'et\'e complexe, on
lui associe deux courants importants en analyse et g\'eom\'etrie
complexe, le courant d'int\'egration sur $\{f=0\}$ et la valeur
principale associ\'ee \`a $f$ (voir [7] et [5]) de la mani\`ere
suivante:
$$<V_p(f),\phi>=\lim_{\varepsilon\rightarrow 0}\int_{\vert f\vert \geq\varepsilon}{\phi\over
f},$$
$$<[f=0],\phi>=\int_{f=0}\phi.$$
Ces d\'efinitions m\'eritent quelques commentaires : pour ce qui
concerne le courant d'int\'egration, d'apr\`es les travaux de
P.Lelong, il faut comprendre que l'int\'egrale est prise sur les
points r\'eguliers de $\{f=0\}$, et que cette derni\`ere est une
int\'egrale g\'en\'eralis\'ee convergente, car les singularit\'es de
$\{f=0\}$ sont de codimension au moins $1$ dans $\{f=0\}$ (voir par
exemple [5] pour les d\'efinitions rigoureuses); l'existence de la
valeur principale est plus d\'elicate et r\'eside essentiellement
sur le th\'eor\`eme de la r\'esolution des singularit\'es (voir [6],
[7]). Il y a un prix \`a payer \`a l'utilisation de ce r\'esultat :
aucune information sur l'ordre de cette distribution. Ces
d\'efinitions peuvent \^etre g\'en\'eralis\'ees \`a la codimension
sup\'erieure, voir [6], [10], pour la valeur principale et [5] pour
le courant d'int\'egration. Dans ce qui suit, nous allons parler des
g\'en\'eralisations de la valeur principale, les courants
r\'esiduels; dans [10], apparaissent les courants suivants, $X_I^J$,
pour $I$ et $J$ deux sous-ensembles disjoints de $\{1,\cdots,p\}$ et
$\theta\rightarrow\varepsilon(\theta)$ un pav\'e admissible (voir
[10] et [6], pour la d\'efinition). Si $\{f_1,\cdots,f_p\}$ est une
intersection compl\`ete,
$$<X_I^J,\phi>=\lim_{\theta\rightarrow 0}\int_{\{\vert
f_I\vert\geq \varepsilon_I(\theta),\vert
f_J\vert=\varepsilon_J(\theta)\}}{\phi\over f_If_J}.$$Ces derniers
courants sont tr\`es importants (qui sont toujours obtenus par la
r\'esolu-tion des singularit\'es) car il est montr\'e dans [10], que
$X_{\emptyset}^{\{1,\cdots,p\}}$ donne la description de l'id\'eal
engendr\'e par $(f_1,\cdots,f_p)$, et que les courants $X_I^J$ donne
une d\'ecomposition "explicite" d'une fonction appartenant \`a
l'id\'eal engendr\'e par $(f_1,\cdots,f_p)$.

Il est bien connu (voir [10]) que la caract\'erisation de l'id\'eal
engendr\'e par $(f_1,\cdots, f_p)$ est de nature cohomologique : si
l'on d\'efinit
$$Res_{f}(\phi):=\int_{\{\vert f_1\vert =\varepsilon_1,\cdots,\vert
f_p\vert =\varepsilon_p\}}{\phi\over f_1\cdots f_p},$$ avec $\phi$
une $(n,n-p)$-forme test $\bar\partial$-ferm\'ee au voisinage de
$f^{-1}(0)$, alors $Res_f$ est ind\'ependant de $\varepsilon$ assez
petit, et nous avons le r\'esultat : \begin{theo}([10])Si $g$ est
holomorphe au voisinage de $z_0\in \{f=0\}$, alors $g\in I_{f}$
(localement au voisinage de $z_0$) si et seulement si $gRes_f=0.$
\end{theo}
Par le th\'eor\`eme d'Hahn-Banach, l'action de $Res_f$ peut \^etre
prolong\'e en un courant $X$ sur les $(n,n-p)$-formes tests qui va
d\'ecrire $I_f$, mais cette approche n'est pas constructive car elle
repose sur l'axiome du choix. Dans [10], on construit une extension
particuli\`ere de $Res_f$, $X_{\emptyset}^{\{1,\cdots,p\}}$, qui est
minimale dans le sens o\`u $\bar gX_{\emptyset}^{\{1,\cdots,p\}}=0$,
si $g\in\sqrt{I_f}$; malheureusement, elle est obtenue par le
th\'eor\`eme de r\'esolution des singularit\'es, ce qui ne rend pas
le proc\'ed\'e plus explicite pour autant.

Dans ce papier, nous proposons la construction d'une extension de
$Res_f$ parfaitement explicite (uniquement avec le th\'eor\`eme de
pr\'eparation de Weierstrass), mais qui n'est pas minimale au sens
pr\'ec\'edent, et donc diff\`ere du courant de Coleff-Herrera.

Dans [8], nous avions construit une extension explicite de $Res_f$
dans le cas de la codimension 1 et dans certains cas particuliers
d'intersection compl\`ete, pour la codimension sup\'erieure. On
v\'erifie ais\'ement - comme nous l'a signal\'e le referee - que le
courant $X$ obtenu pour, par exemple $f=(z_1^2-z_2^2)$, n'a pas la
propri\'et\'e d'extension minimale et n'est donc pas le courant
r\'esiduel classique $\bar\partial V_p(f)$.

Evidemment plus les courants pr\'ec\'edents sont obtenus de
mani\`ere explicite, plus on obtient d'informations sur l'id\'eal et
sur la d\'ecomposition d'un \'el\'ement particulier; c'est pourquoi,
le th\'eor\`eme de r\'esolution des singularit\'es peut s'av\'erer
un obstacle difficilement surmontable par exemple, pour obtenir des
d\'ecompositions dans certains espaces de r\'egularit\'es pour les
fonctions holomorphes (par des m\'ethodes $L^2$, H.Skoda a obtenu un
th\'eor\`eme de ce type important dans les espaces $L^2$ \`a poids).

Ici, nous proposons une autre approche plus \'el\'ementaire, pour
construire des courants ayant les m\^emes propri\'et\'es de ceux de
Coleff-Herrera-Passare (en ce qui concerne la description de $I_f$
et de la d\'ecomposition dans $I_f$) dans le cas d'une intersection
compl\`ete quelconque. Cette approche repose essentiellement sur le
th\'eor\`eme de pr\'eparation de Weierstrass et donne donc l'ordre
des courants construits. Cette construction est la suite d'articles
pr\'ec\'edents (voir [8], [9]) qui donnaient ces courants dans des
intersections compl\`etes particuli\`eres. Pour rendre l'article
plus lisible, nous commencons la construction pour une intersection
compl\`ete de codimension 2, le cas g\'en\'eral est obtenu avec la
m\^eme philosophie dans les parties suivantes.
\bigskip

\centerline{Liste des notations.}
\bigskip

-Pour $\phi=\sum\phi_{IJ}dz_I\wedge d\bar z_J$ une $(p,q)$-forme
diff\'erentielle, on note $\bar\partial\phi$ la partie de bidegr\'es
$(p,q+1)$ de sa diff\'erentelle ext\'erieure.
\medskip

-On notera $\bar\partial_l(\phi)$, la forme diff\'erentielle de
bidegr\'es $(p,q)$ d\'efinie par l'expression:
$$\bar\partial_l(\phi)=\sum{\partial\phi_{IJ}\over\partial \bar z^l}dz_I\wedge d\bar
z_J.$$

\medskip

-Pour $f$ une fonction holomorphe v\'erifiant les conditions du
th\'eor\`eme de Weierstrass par rapport \`a la variable $z_l$, on
note $P_l^f$ son polyn\^ ome de Weierstrass par rapport \`a $z_l$
associ\'e. En g\'en\'eral, nous d\'esignerons par $N_l$ le degr\'es
de $P_l^f$ par rapport \`a $z_l$.

\medskip

-Pour $(f_1,\cdots,f_p)$ une intersection compl\`ete au voisinage de
$0$, nous noterons \'egalement $(f)$.

\medskip

-Si $P_l,Q_l$ sont deux polyn\^omes de Weierstrass par rapport \`a
la variable $z_l$, $R_l(P_l,Q_l)$ d\'esigne le r\'esultant de $P_l$
et $Q_l$ par rapport \`a la variable $z_l$ (il est \`a noter que
$R_l$ ne d\'epend pas de la variable $z_l$).
 \bigskip

 Dans [8] et [9], nous avions d\'egag\'e "les propri\'et\'es fonctionnelles" que devaient v\'erifier
 une famille de courants attach\'ee \`a $(f_1,f_2)$ afin de
 d\'ecrire $I_f$. Plus pr\'ecis\'ement, rappelons le th\'eor\` eme de [8] :
\begin{theo}\label{thm:main}Soit $(\theta_1,\theta_2)$ une intersection compl\` ete d\'efinie au
voisinage de $0$ ; s'il existe une famille de courant $(X_1,X_2)$
v\'erifiant les propri\'et\'es ci-dessous :
$$\theta_1X_1=1,\ \theta_2X_2=\bar\partial X_1,\ \theta_1X_2=0,$$
alors $g\in I_{loc}(\theta_1,\theta_2)$ si et seulement si $g\bar\partial X_2=0$.
\end{theo}
\begin{rem}Le th\'eor\`eme pr\'ec\`edent (et le th\'eor\`eme 4.1 par la suite,
pour la codimension sup\'erieure) est encore vrai d\`es que l'on
peut trouver $(X_1,X_2)$ v\'erifiant les deux propri\'et\'es
ci-dessus, sans que $(\theta_1,\theta_2)$ soit forc\'ement une
intersection compl\`ete - la preuve de [8] repose uniquement sur les
propri\'et\'es fonctionnelles des courants -. Nous le mentionnons
dans le cas d'une intersection compl\`ete car il s'agit du cadre
naturel de son application : en effet dans le cas non intersection
compl\`ete, il est clairement impossible, en g\'en\'eral, de
produire des courants avec ces propri\'et\'es - pour s'en
convaincre, consid\'erer la situation  $f=(z_1,z_1)$ -. De m\^eme,
le th\'eor\`eme de d\'ecomposition (3.1) est valable - par la m\^eme
preuve produite dans [8] - sans l'hypoth\`ese d'intersection
compl\`ete.
\medskip

On nous a fait remarquer que l'on peut obtenir les th\'eor\`emes 1.2
et 4.1 en utilisant la th\'eorie g\'en\'erale d\'evelopp\'ee dans
[1]et [2]. Mais \`a notre avis, l'approche de [8] est plus
\'el\'ementaire et suffit \`a notre construction.
\medskip

Enfin, nous tenons \`a signaler que dans le cas d'une intersection
compl\`ete, le calcul sur les courants r\'esiduels classiques (les
courants $X_I^J$ d\'efinis pr\'ec\'edemment) d\'evelopp\'e dans
[11], [4] et [13], assure qu'ils v\'erifient les conditions du
th\'eor\`eme 1.2, et en codimension sup\'erieure les conditions (1)
et (2) du th\'eor\`eme 4.1.
\end{rem}

Ici, nous voulons produire d'autres courants (que les courants de
Coleff-Herrera-Passare) v\'erifiant les hypoth\`eses du th\'eor\`eme
1.2 - sans utiliser le th\'eor\`eme d'Hironaka -, et de la mani\`ere
la plus constructive possible. En g\'en\'eral, il est tr\`es
difficile de construire directement de tels courants, sauf dans
certains cas particuliers d'intersections compl\`etes (voir [8] et
[9]). L'id\'ee est donc de toujours  se ramener \`a cette situation
(voir, section : Construction d'une intersection compl\` ete
adapt\'ee et ce qui suit pour la codimension 2).
\medskip

Soit $(f_1,f_2)$ une intersection compl\` ete d\'efinie au voisinage
de $0$ dans $\C^n$ avec $0\in\{f_1=f_2=0\}$ ; quitte \` a faire un
changement lin\'eaire de coordonn\'ees, on peut supposer que
$(f_1,f_2)$ v\'erifient les conditions du th\'eor\` eme de
pr\'eparation de Weierstrass par rapport \` a la variable $z_1$, et
obtenir ainsi : $f_1=U_{f_1}P^{f_1}_1$ et $f_2=U_{f_2}P^{f_2}_1$ au
voisinage de $0$, avec $U_i$ des unit\'es.
\medskip

Nous avons donc :
$$R_1(P^{f_1}_1,P^{f_2}_1)(Z^{'})=a_1(z)f_1(z)+a_2(z)f_2(z)\ (*),$$
avec $a_i$ des fonctions holomorphes au voisinage de $0$ ; on peut
remarquer que si $(f_1,f_2)$ est une intersection compl\` ete alors
$(f_1,R_1)$ l'est \'egalement. Clairement $g\in I_{loc}(f_1,f_2)$ si
et seulement si $gdet(A)\in I_{loc}(f_1,R_1)$, o\`u $A$ est la
matrice holomorphe de passage de $(f_1,f_2)$ \`a $(f_1,R_1)$ ; on
peut \'ecrire, en utilisant $A$, $(*)$ sous la forme :
 $$ R_1(Z^{'})   =a_1(z)f_1(z)+det(A)f_2(z)\ (**).$$

Ainsi, nous allons appliquer le th\'eor\`eme 1.2, non pas
directement \` a $(f_1,f_2)$, mais \`a l'intersection compl\` ete
$(f_1,R_1)$ qui elle est adapt\'ee.

\section{Description de $I_{loc}(f_1,f_2)$.}

Commencons par nous placer - quitte \`a faire \` a nouveau un
changement lin\'eaire de coordonn\'ees, mais uniquement sur les
variables $(Z^{'})$, cette fois - dans un syst\`eme de variables,
pour lequel $R_1$ v\'erifie les conditions du th\'eor\`eme de
Weierstrass par rapport \` a la variable $z_2$ (ceci nous sera utile
par la suite).

\medskip

D\'efinissons le courant $X_{1}^{\Gamma}$, pour $\Gamma\in \N$, par
son action sur  les $(n,n)$-formes tests au voisinage de $0$ :
$$<X_{1}^{\Gamma},\phi>=\int_{\C^n}{(\overline{ P_{1}^{f_1}})^{\Gamma}\over
f_1}\bar\partial_1^{\Gamma N_1}(\phi),$$ o\` u $N_1$ est le degr\'es
de $P_{1}^{f_1}$. Il est alors assez facile de voir ([8]), que
$X_{1}^{\Gamma}$ v\'erifie l'\'equation - modulo une constante de
normalisation - :
$$f_1 X_{1}^{\Gamma}=1,\ \forall\ \Gamma\in \N;$$
de m\^eme, il est ais\'e d'obtenir l'expression de  $\bar\partial
X_{1}^{\Gamma}$ sur les $(n,n-1)$-formes :
$$<\bar\partial X_{1}^{\Gamma},\phi>=\int_{\C^n}{\bar\partial(\overline{ P_{1}^{f_1}})^{\Gamma}\over
f_1}\bar\partial_1^{\Gamma N_1}(\phi).$$

\medskip

Maintenant, construisons $X_{2}^{\Gamma}$, un $(n,n-1)$-courant, par
la formule :
$$<X_{2}^{\Gamma},\phi>=\int_{\C^n}{\overline{P_{2}^{R_1}}\over
R_1}\bar\partial_{2}^{N_2}\bigg({\bar
\partial(\overline{ P_{1}^{f_1}})^{\Gamma}\over f_1}\bar\partial_1^{\Gamma
N_1}(\phi)\bigg),$$ o\` u $N_2$ est le degr\'es du polyn\^ome de
Weierstrass associ\'e \`a $R_1$ - $P_{2}^{R_1}$ -, et enfin $\Gamma$
est choisi de sorte que l'int\'egrande soit une fonction r\'eguli\`
ere.

\medskip
\centerline{Fait 1 : $R_1X_{2}^{\Gamma}=\bar\partial
X_{1}^{\Gamma}$}
\medskip

En effet, explicitons l'action du courant \` a droite de
l\'egalit\'e :
$$<R_1X_{2}^{\Gamma},\phi>=\int_{\C^n}\overline{P_{2}^{R_1}}
\bar\partial_{2}^{N_2}\bigg({\bar
\partial(\overline{ P_{1}^{f_1}})^{\Gamma}\bar\partial_1^{\Gamma
N_1}(\phi)\over f_1}\bigg),$$ et donc, \` a l'aide de
$N_2$-int\'egration par partie - par rapport \` a la variable $z_2$
-, on obtient l'\'egalit\'e.

\medskip \centerline{Fait 2 :
$f_1X_{2}^{\Gamma}=0$}
\medskip

En effet, nous avons les \'egalit\'es :
\begin{align*}
 <f_1X_{2}^{\Gamma},\phi> & =\int_{\C^n}{\overline{P_{2}^{R_1}}\over
R_1}\bar\partial_{2}^{N_2}\bigg(\bar
\partial(\overline{ P_{1}^{f_1}})^{\Gamma}\bar\partial_1^{\Gamma
N_1}(\phi)\bigg)\\
& =\sum_{\alpha+\beta=N_2}\int_{\C^n}{\overline{P_{2}^{R_1}}\over
R_1}\bar\partial_{2}^{\alpha}\bigg(\bar
\partial(\overline{ P_{1}^{f_1}})^{\Gamma}\bigg)\bar\partial_{2}^{\beta}\bigg(\bar\partial_1^{\Gamma
N_1}(\phi)\bigg)\\
& =\sum_{\alpha+\beta=N_2}\int_{\C^n}{\overline{P_{2}^{R_1}}\over
R_1}\bar\partial_{2}^{\alpha}\bigg(\bar
\partial(\overline{
P_{1}^{f_1}})^{\Gamma}\bigg)\bar\partial_{1}^{\Gamma
N_1}\bigg(\bar\partial_2^{\beta}(\phi)\bigg);
\end{align*}
le facteur ${\overline{P_{2}^{R_1}}\over R_1}$ ne d\'epend pas de la
variable $z_1$, ceci entra\^\i ne - en int\'egrant $\Gamma N_1$-fois
par partie par rapport \` a la variable $z_1$, et car
$\bar\partial_{2}^{\alpha}\bigg(\bar
\partial(\overline{
P_{1}^{f_1}})^{\Gamma}\bigg)$ est de degr\'es strictement
inf\'erieur \` a $\Gamma N_1$ - l'\'egalit\'e souhait\'ee.

\medskip

Nous pouvons maintenant d\'ecrire l'id\'eal local, en $0$,
engendr\'e par $(f_1,f_2)$:
\begin{theo}\label{thm:main}$g\in I_{loc}(f_1,f_2)$ si et seulement
si le courant $gdet(A)\bar\partial(X_{2}^{\Gamma})=0$.
\end{theo}
Remarque :

Le proc\'ed\'e pour obtenir le courant $X_{2}^{\Gamma}$ est
parfaitement explicite-la matrice $A$ et le r\'esultant $R_1$
s'obtiennent en d\'eroulant l'algorithme d'Euclide pour les
polyn\^omes de Weierstrass $(P^{f_1}_1,P^{f_2}_1)$ -.

\medskip

Preuve du th\'eor\`eme :

comme nous avons vu lors de la section 1, $g\in I_{loc}(f_1,f_2)$ si
et seulement si $gdet(A)\in I_{loc}(f_1,R_1)$ ; mais nous sommes
capables de d\'ecrire cet id\'eal-car la famille
$(X_{1}^{\Gamma},X_{2}^{\Gamma})$ v\'erifie les conditions du
th\'eor\` eme 1.2-, et donc $g\in I_{loc}(f_1,f_2)$ si et seulement
si $gdet(A)\bar\partial (X_{2}^{\Gamma})=g\bar\partial
(det(A)X_{2}^{\Gamma})=0$.

\section{Formule de d\'ecomposition dans $I_{loc}(f_1,f_2)$. }

Dans cette section, nous supposerons que
$\theta=(\theta_1,\theta_2)$ est d\'efinie au voisinage de $\bar
{B_r}$, la boule ferm\'ee de $\C^n$ de centre $0$ et rayon $r$.

Nous consid\`erons \'egalement des fonctions $\theta^i_1(\zeta,z)$
et $\theta^i_2(\zeta,z)$ holomorphes sur $B_r\times B_r$ telles que
: $$ \theta_1(z)-\theta_1(\zeta)=\sum_i
\theta^i_1(\zeta,z)(z_i-\zeta_i)\ \
\theta_2(z)-\theta_2(\zeta)=\sum_i
\theta^i_2(\zeta,z)(z_i-\zeta_i),$$ autrement dit une
d\'ecomposition de Hefer de $\theta_1$ et $\theta_2$ sur $B_r$. Si
$(X_1,X_2)$ est une famille de courants v\'erifiant les conditions
du th\'eor\`eme 1.2, nous avons le r\'esultat plus pr\'ecis de
d\'ecomposition dans $I_{\theta_1,\theta_2}$ (voir [8]).
\begin{theo}\label{thm:main} Soit $g$
holomorphe sur $\bar {B_r}$ avec $g\bar\partial X_2=0$. Alors, il
existe $P_1(\zeta,z)$, $P_2(\zeta,z)$ deux noyaux int\'egraux
holomorphes en $z$, ne d\'ependant que de $B_r$, tels que $\forall
z\in B_r$ : \begin{align*} g(z) =\theta_1(z)<gX_1,P_1(.,z)> &
+\theta_1(z)<\sum_i \theta^i_2(\zeta,z)d\zeta_i\wedge
X_2,P_2(.,z)>\\
& +\theta_2(z)<\sum_i \theta^i_1(\zeta,z)d\zeta_i\wedge
X_2,P_2(.,z)>.\end{align*}
\end{theo}
Consid\'erons le courant $det (A)X_{2}^{\Gamma}$ construit \` a la section
pr\'ec\'edente, alors nous avons d'une part :
 $$ f_1det(A)X_{2}^{\Gamma}=0,$$
et d'autre part,en utilisant les propri\'et\'es des courants
$(X_{1}^{\Gamma},X_{2}^{\Gamma})$,
$$f_2det(A)X_{2}^{\Gamma}
=(R_1-a_1f_1)X_{2}^{\Gamma}=R_1X_{2}^{\Gamma}=\bar\partial
X_{1}^{\Gamma}.$$ Par cons\'equent, la famille de courants
$(X_{1}^{\Gamma},det(A)X_{2}^{\Gamma})$ v\'erifie, pour $(f_1,f_2)$,
les conditions d'application du th\'eor\` eme 3.1 - \`a des
constantes de normalisation pr\`es -, et l'on obtient une formule de
d\'ecomposition pour $g\in I_{loc}(f_1,f_2)$.

\section{Cas de la codimension sup\'erieure \` a $1$.}

Enoncons les deux r\'esultats de [8], sur lesquels repose cette
construction.

Soit $(\theta_1,\cdots,\theta_p)$ une intersection compl\` ete au
voisinage de $0$ ; supposons qu'il existe une famille de courants
$(X_1,\cdots,X_p)$ v\'erifiant les propri\'et\'es $(1)$ et $(2)$
suivantes:
$$\begin{cases} & \ \ (1): \theta_1X_1=1,\ \ \theta_jX_j=\bar\partial X_{j-1},\ \ \forall
j\in\{2,\cdots,p\},\\
& \ \ (2): \forall i\in\{2,\cdots,p\},\ \ \theta_jX_i=0,\ \ \forall
j\in\{1,\cdots,i-1\}.\end{cases}$$

Alors, nous avons la description de
$I_{loc}(\theta_1,\cdots,\theta_p)$ : \begin{theo}\label{thm:main}
Sous les conditions pr\'ec\'edentes :
$$g\in I_{loc}(\theta_1,\cdots,\theta_p)\ \ \hbox{si et seulement si}\ \
g\bar\partial X_p=0.$$
\end{theo}
\begin{rem}Comme pour le th\'eor\`eme 1.2, le r\'esultat est valable
d\`es que l'on peut trouver des courants v\'erifiant (1) et (2),
sans l'hypoth\`ese d' intersection compl\`ete - la preuve est la
m\^eme que celle de [8]-.
\medskip

De m\^eme, le th\'eor\`eme suivant de d\'ecomposition est valable
d\`es que les courants v\'erifient (1) et (2) avec la m\^eme preuve
de [8].
\end{rem}

 En fait, nous avons un r\'esultat  de d\'ecomposition
dans $I_{loc}(\theta_1,\cdots,\theta_p)$, comme dans la cas de la
codimension $2$ :

\begin{theo}\label{thm:main} Soit $(\theta_1,\cdots,\theta_p)$, une intersection compl\` ete sur
une boule de centre $0$ et de rayon $r$, $B_r$, soit
$(X_1,\cdots,X_p)$ une famille de courant v\'erifiant les conditions
$(1)$ et $(2)$ ci-dessus, et enfin, $g$ holomorphe sur $\bar B_r$
avec $g\bar\partial X_p=0$. Alors, il existe des noyaux int\'egraux,
$P_1(\zeta,z),\cdots,P_p(\zeta,z)$, holomorphes en $z$, ne
d\'ependant que de $B_r$, tels que, $\forall\ z\in B_r$ :

$$g(z)=\sum_{i\geq j}C_{ij}
\theta_j(z)<g(\zeta)B_{i}^{j}(\zeta,z)\wedge X_i,P_i(.,z)>,$$
avec les notations suivantes :
$$\begin{cases} & B_j(\zeta,z)=b_1(\zeta,z)\wedge\cdots\wedge
b_j(\zeta,z),\\
&
B_{j}^{l}=b_1(\zeta,z)\wedge\cdots\wedge\widehat{b_l(\zeta,z)}\wedge\cdots\wedge
b_j(\zeta,z),\ \l\leq j\ \ \hbox{et}\ B_{1}^{1}(\zeta,z)=1,\\
& b_l(\zeta,z)=\sum_{i=1}^{n}b_{l}^{i}(\zeta,z)d\zeta_i\ \
\hbox{o\`u}\ \ \theta_l(z)-\theta_l(\zeta)=
\sum_{i=1}^{n}b_{l}^{i}(\zeta,z)(z_i-\zeta_i).\end{cases}$$
\end{theo}
Si $(f_1,\cdots,f_p)$ est une intersection compl\`ete donn\'ee au
voisinage de $0$, comme dans la section pr\'ec\'edente, nous
n'allons pas construire directement-car cela est \`a priori pas
\'evident...-la famille $(X_1,\cdots,X_p)$ associ\'ee ; mais, nous
allons passer par une intersection compl\`ete interm\'ediaire,
construite \` a partir de la premi\`ere, pour laquelle, il est
facile de construire la famille de courant adhoc.

\section{Construction d'une intersection compl\` ete adapt\'ee.}

Nous allons proc\'eder par induction. Soit $U_1$ le vecteur de
composantes $(f_i)_{1\leq i\leq p}$; pour construire $U_2$,
commencons par faire un changement de coordonn\'ees lin\'eaires,
afin que les $f_i$ v\'erifient le th\'eor\` eme de pr\'eparation de
Weierstrass par rapport \` a la variable $z_1$. Dans ces nouvelles
coordonn\'ees, $U_2$ est le vecteur de composantes $$(f_1,
R_{1}(P_{1}^{f_1},\sum_{l=1}^{2}\lambda_{2,2}^{l}P_{1}^{f_l}),\cdots,
 R_{1}(P_{1}^{f_1},\sum_{l=1}^{j}\lambda_{2,j}^{l}P_{1}^{f_l}),\cdots,R_{1}(P_{1}^{f_1},
 \sum_{l=1}^{p}\lambda_{2,p}^{l}P_{1}^{f_l})),$$
 o\`u $j$ est l'indice qui indique la position dans le p-uplet $U_2$, les $\lambda_{2,j}^{l}$
 sont des scalaires choisis de mani\`ere \`a ce que $U_2$ soit une intersection compl\`ete en
 $0$; ceci est possible car $(f_1,\cdots,f_p)$ l'est.\footnote{Il faut
remarquer que $(f_1, R_{1}(P_{1}^{f_1},P_{1}^{f_2}),\cdots,
 R_{1}(P_{1}^{f_1},P_{1}^{f_p}))$ n'est pas une intersection
 compl-\`ete en g\'en\'eral-m\^eme si $(f_1,\cdots,f_p)$ l'est-; il est donc
n\'ecessaire de choisir des scalaires, $(\lambda_{2,j}^l)(2\leq
j\leq p, 1\leq l\leq j),$
 pour que
 $$(f_1,R_{1}(P_{1}^{f_1}, \lambda^{l}_{2,2}P_{1}^{f_2}),\cdots,R_{1}(P_{1}^{f_1},
\lambda_{2,j}^lP_{1}^{f_l}),\cdots, R_{1}(P_{1}^{f_1},
\lambda_{2,p}^lP_{1}^{f_l})),$$ o\`u l'on somme sur l'indice du haut
variant de $1$ \`a l'indice du bas \`a droite, soit elle, une
intersection compl\`ete. En voici un exemple :
$$(f_1(z)=(z_1+z_2)(z_1+z_3+z_2z_3),f_2(z)=z_1+z_2^2,f_3(z)=z_1+z_2+z_3);$$
un calcul \'el\'ementaire conduit aux expressions suivantes :
$$R(f_1,f_2)=z_2(z_2^3-z_2^2-z_2^2z_3+z_3),
R(f_1,f_3)=z_2z_3(1-z_3).$$ Ainsi sur $\{z_2=0\}$, $R(f_1,f_2)$ et
$R(f_1,f_3)$ sont tous deux nuls...mais en prenant une perturbation
lin\'eaire de $f_3$, soit en posant $g=af_2+bf_3$ avec $a$, $b$ des
constantes complexes, on obtient $(f_1, R(f_1,f_2), R(f_1,g))$ qui
elle, est une intersection compl\`ete. }
\bigskip

En effet, soit $(g_1,\cdots,g_p)$ des polyn\^omes de Weierstrass par
rapport \`a la variable $z_1$, d\'efinissant une intersection
compl\`ete en z\'ero : il suffit de montrer que si
$(g_1,R_1(g_1,g_2),\cdots,R_1(g_1,g_i))$ est une intersection
compl\`ete au voisinage de z\'ero et
$(g_1,R_1(g_1,g_2),\cdots,R_1(g_1,g_i),R_1(g_1,g_{i+1}))$ non, alors
il existe $(\lambda_l)$ tels que
$$(g_1,R_1(g_1,g_2),\cdots,R_1(g_1,g_i),R_1(g_1,\sum_{l=1}^{i+1}\lambda_lg_{l}))$$
soit de dimension $n-i-1$. Notons pour $l\leq p$
\begin{align*}
&
W_l:=\{R_1(g_1,g_2)=\cdots=R_1(g_1,g_l)=0\}\subset\C^{n-1}=(z_2\cdots,z_n)\
,\\
& V_l:=\{g_1=\cdots=g_l=0\};\end{align*} $W_{i+1}$ ne d\'epend pas
de la variable $z_1$ donc - si $W_{i+1}\cap\{g_{1}=0\}$ n'est pas
une intersection compl\`ete en z\'ero -, alors $R_1(g_1,g_{i+1})$
est identiquement nulle sur une composante irr\'eductible $C^{j_0}$
de $W_i$. Consid\'erons la projection $\pi_1$ :
$$\{g_1=0\}\rightarrow \C^{n-1}\
,(z_1,\cdots,z_n)\rightarrow(z_2,\cdots,z_n);$$ nous avons
$\pi_1(V_i)\subset W_i$, $V_i\subset \pi_1^{-1}(W_i)$ et de plus,
$dim(V_i)=dim(\pi_1^{-1}(W_i))=n-i$, ce qui entra\^\i ne que
$\pi_1^{-1}(W_i)$ poss\`ede quelques composantes irr\'eductibles de
plus que $V_i$ :
$$\pi_1^{-1}(W_i)=(\cup_k\Gamma_
{V_i}^k)\cup(\cup_j\Gamma^j),$$ o\`u $\Gamma_ {V_i}^l$ sont les
composantes irr\'eductibles de $V_i$ et $\Gamma^j$ les autres!
\bigskip

Maintenant $R_1(g_1,g_{i+1})$ identiquement nulle sur $C^{j_0}$
implique $g_{i+1}$ identiquement nulle sur une composante
irr\'eductible de $\pi_1^{-1}(W_i)$; $V_{i+1}$ est une intersection
compl\`ete et donc, $\forall \ k\ dim(\{g_{i+1}=0\}\cap
\Gamma^k_{V_i})=n-i-1$, ce qui entra\^\i ne l'existence de $j_1$ tel
que $g_{i+1}=0$ sur $\Gamma^{j_1}$, identiquement.
\bigskip

Consid\`erons l'ensemble $A$:
$$A:=\{j/g_{i+1}=0\ sur\ \Gamma^j\};$$
pour tout $j\in A$, $dim(\Gamma^j\cap V_i)=n-i-1$, il existe donc
$(\lambda_l)$ - $\lambda_{i+1}$ peut \^etre choisi non nul - tels
que $\sum_{l=1}^{i+1}\lambda_lg_l$ est non identiquement nulle sur
$\Gamma^j$ pour tout $j\in A$ et sur $\Gamma^j$ pour tout $j$.
Clairement  une  combinaison lin\'eaire de ce type ne peut \^etre
identiquement nulle sur un $\Gamma^k_{V_i}$, on a donc
$\sum_{l=1}^{i+1}\lambda_lg_l$ non identiquement nulle sur une
composante irr\'eductible de $\pi_1^{-1}(W_i)$, ainsi
$$dim(\{R_1(g_1,\sum_{l=1}^{i+1}\lambda_lg_l)=0\}\cap W_i)=n-i-1.$$
ce qui termine la preuve de l'assertion.
\bigskip

On note en outre que les composantes $U_2^{j}$ de $U_2$, $j>1$, sont
des fonctions holomorphes des variables $(z_2,\cdots,z_p)$.
\bigskip

Supposons $U_{i-1}$ un $p$-vecteur donn\'e, d\'efinissant une
intersection compl\` ete au voisinage de $0$, de coordonn\'ees
$U_{i-1}^{j}$ avec $U_{i-1}^{j}$ des fonctions holomorphes des
variables $(z_j,\cdots,z_n)$, v\'erifiant le th\'eor\`eme de
pr\'eparation de Weierstrass par rapport \`a la variable $z_j$ pour
$j<i-1$, et $U_{i-1}^j$ des fonctions holomorphes des variables
$(z_{i-1},\cdots,z_p)$, pour $j\geq i-1$.

Le vecteur $U_i$ est alors obtenu de la mani\` ere suivante-en
effectuant un changement de coordonn\'ees lin\'eaires en les
variables $(z_{i-1},\cdots,z_p)$, de facon \`a ce que les
$U_{i-1}^{j}$, $j\geq i-1$, v\'erifient le th\'eor\` eme de
pr\'eparation de Weierstrass par rapport \` a la variable
$z_{i-1}$-; on choisit des scalaires $\lambda_{i,j}^{l},\ (j\geq
i),\ (l\leq j)$ comme ci-dessus et l'on pose:
 $$U_i=(U_{i-1}^{1},\cdots,U_{i-1}^{i-1},R_{i-1}(P_{i-1}^{U_{i-1}^{i-1}},\sum_{l=1}^{i}\lambda_{i,i}^{l}P_{i-1}^{U_{i-1}^{l}}),
 \cdots, R_{i-1}(P_{i-1}^{U_{i-1}^{i-1}},\sum_{l=1}^{p}\lambda_{i,p}^{l}P_{i-1}^{U_{i-1}^{l}})).$$
\medskip

-On note d'une part, que $U_{i}^{j}$ sont des fonctions holomorphes
en les variables $(z_j,\cdots, z_n)$ pour $j<i$ qui v\'erifient le
th\'eor\`eme de pr\'eparation de Weierstrass par rapport \`a $z_j$,
et $U_{i}^{j}$ sont des fonctions holomorphes des variables
$(z_i,\cdots z_n)$, pour $j\geq i$.

-D'autre part, que $U_i$ d\'efinit une intersection compl\`ete au
voisinage de $0$ si $U_{i-1}$ en d\'efinit une, pourvu que les
$\lambda_{i,j}^{l}$ soient correctement choisis.

\medskip

-Pour finir, on pose $(f_1,R_2,\cdots,R_p):= U_{p}$ et nous avons la
situation suivante : les fonctions $R_i$, sont des fonctions
holomorphes en les variables $(z_i,\cdots z_n)$ v\'erifiant le
th\'eor\`eme de pr\'eparation de Weierstrass par rapport \`a $z_i$.

\section{Construction de la famille $(X_1,\cdots,X_p)$ associ\'ee
\`a $(f_1,R_2,\cdots,R_p)$ }

\'Etant donn\'e une intersection compl\` ete, il n'est pas simple en
g\'en\'eral de construire une famille de courants v\'erifiant les
conditions de la section 4. Par contre, pour une intersection
compl\` ete de la forme $(\theta_1,\cdots,\theta_p)$ avec $\theta_i$
des fonctions holomorphes des variables $(z_i,\cdots,z_n)$
v\'erifiant le th\'eor\`eme de pr\'eparation par rapport \`a $z_i$,
cela est plus ais\'e. Choisissons $(\Gamma_1,\cdots,\Gamma_p)$ une
famille de $p$ entiers naturels de mani\`ere \`a ce que les
fonctions intervenant dans les int\'egrandes soient r\'eguli\`eres,
et posons,

$$\begin{cases}<X_1,\phi> =\int {(\bar
P_{1}^{\theta_1})^{\Gamma_1}\over
\theta_1} & \bar\partial_1^{\Gamma_1N_1}(\phi),\\
<X_i,\phi>  = \int {(\bar p_{i}^{\theta_i})^{\Gamma_i}\over
\theta_i} &  \bar\partial_{i}^{\Gamma_iN_{i}}\big({\bar\partial(\bar
P_{i-1}^{\theta_{i-1}})^{\Gamma_{i-1}}\over
\theta_{i-1}}\bar\partial_{i-1}^{\Gamma_{i-1}N_{i-1}}  \big(\cdots
\\
&  \cdots\big(\bar\partial_{2}^{\Gamma_2 N_2}\big({\bar\partial(\bar
P_1^{\theta_1})^{\Gamma_1}\over
\theta_1}\bar\partial_1^{\Gamma_1N_1}(\phi)\big)\big)\cdots\big).\end{cases}$$
Alors, par des arguments similaires au cas de la codimension $2$, on
v\'erifie sans peine que la famille de courants $(X_1,\cdots,X_p)$
poss\` ede toutes les propri\'et\'es requises.

\bigskip

Soit $A=(a_{i,j})$-$i$ d\'esigne la colonne et $j$ la ligne-la
matrice triangulaire sup\'erieure, dont le pivot $a_{1,1}$ vaut $1$,
qui transforme $(f_1,\cdots,f_p)$ en $(f_1,R_2,\cdots,R_{p})$; soit
$\delta_l:=\prod_{i=1}^{l}a_{i,i}$, pour $l\in \{1,\cdots,p\}$ et
enfin, $(X_1,\cdots,X_p)$ les courants v\'erifiant les hypoth\` eses
du th\'eor\` eme 4.1 pour l'intersection compl\` ete
$(f_1,R_2,\cdots, R_{p})$.
\begin{cor}Soit $(f_1,\cdots,f_p)$ une
intersection compl\` ete. Alors, nous avons $g\in
I_{loc}(f_1,\cdots,f_p)$ si et seulement si $gdet(A)\bar\partial
X_p=0$.
\end{cor}

Preuve du corollaire :

\medskip

Nous avons $g\in I_{loc}(f_1,\cdots,f_p)$ si et seulement si
$gdetA\in I_{loc}(f_1,R_2,\cdots,R_p)$-on peut par exemple utiliser
la formule de transformation des courants r\'esiduels classiques
(voir [3])\footnote{Remarquons que la loi de transformation des
courants classiques de Coleff-Herrera ([3]) assure que, si $(f)$,
$(g)$ sont deux intersections compl\`etes et $(g)=A(f)$ avec $A$ une
$p\times p$-matrice holomorphe, alors
$$\bar\partial\big[{1\over
f_1}\big]\wedge\cdots\wedge\bar\partial\big[{1\over
f_p}\big]=det(A)\bar\partial\big[{1\over
g_1}\big]\wedge\cdots\wedge\bar\partial\big[{1\over g_p}\big].$$
C'est cette propri\'et\'e qui est particuli\`ere aux courants
r\'esiduels classiques dans l'ensemble des courants d\'ecrivant
l'id\'eal local de $f$. Il serait int\'eressant de voir si cette
propri\'et\'e les caract\'erise totalement dans cet ensemble...}-.
Pour conclure, il suffit d'appliquer le th\'eor\` eme 4.1 \` a
$(f_1,R_2,\cdots, R_p)$.

\bigskip

Comme \`a la section pr\'ec\'edente, nous allons obtenir une
d\'ecomposition dans l'id\'eal local \`a partir de la famille
$(X_1,\cdots,X_p)$.

Nous avons $\sum_{i}a_{i,p}f_i=R_p$, et donc
$$\begin{cases} & \delta_{p}f_p
=\delta_{p-1}R_p-\delta_{p-1}\sum_{i=1}^{p-1}a_{i,p}f_i,\\
&\delta_{p}f_pX_p=\delta_{p-1}\bar\partial
X_{p-1}-\sum_{i=1}^{p-1}a_{i,p}\delta_{p-1}f_iX_p.\end{cases}$$
Comme $i\in\{1,\cdots,p-1\}$ et que $A$ est triangulaire, alors
$\sum a_{i,p}\delta_{p-1}f_i$ s'exprime en fonction de
$f_1,R_2,\cdots,R_{p-1}$, et donc
$$\delta_{p}f_pX_p=\delta_{p-1}\bar\partial X_{p-1}.$$ De plus, pour tout
$l<p$-toujours car $A$ est triangulaire-
$$f_l\delta_p=f_{l}a_{l,l}\delta^{l}_{p}=\delta_{p}^{l}R_l-\sum_{i=1}^{l-1}a_{i,l}\delta_{p}^{l}f_i,$$
o\`u $\delta_{p}^{l}$ signifie que l'on omet $a_{l,l}$ dans le
produit. Par cons\'equent, $f_l\delta_p$ s'exprime en fonction de
$f_1,R_2,\cdots,R_l$, ce qui entra\^\i ne $f_l\delta_pX_p=0$, pour
tout $l<p$. On v\'erifie de la m\^eme mani\`ere que
$f_{p-1}\delta_{p-1}X_{p-1}=\delta_{p-2}\bar\partial X_{p-2}$ et
$f_l\delta_{p-1}X_{p-1}=0,$ pour tout $l<p-1$. Au bout du compte, il
vient sans difficult\'e que la famille $(X_1,
\delta_{2}X_2,\cdots,\delta_{p}X_p)$ v\'erifie toutes les conditions
requises du th\'eor\` eme 4.2, et nous avons donc une
d\'ecomposition, en termes de courants \'el\'ementaires, dans
$I_{loc}(f)$.

\section{Quelques exemples.}
\bigskip
Consid\'erons $(f_1,f_2)$, l'intersection compl\` ete au voisinage de $0$
d\'efinie par :
$$(z_1^{2}+z_2z_3,z_1^{2}+z_3z_1+z_2^{2}-z_3^{3});$$
des  calculs \'el\'ementaires donnent les expressions suivantes :
$$R_1(f_1,f_2)=z_2z_3^{3}+(z_2^2-z_3^3-z_2z_3)^2.$$
En posant $P(z)=z_1-z_2^2+z_3^3+z_2z_3$, nous avons :
$$R_1(f_1,f_2)=(z_3^2+P(z))f_1(z)-P(z)f_2(z):=
h_1(z)f_1(z)+h_2(z)f_2(z).$$ En appliquant les r\'esultats des
sections 1 et 2, nous obtenons la famille $(X_1,h_2X_2)$ qui
d\'ecrit l'id\'eal associ\'e \` a $(f_1,f_2)$ :
$$\begin{cases} & <X_1,\phi>=\int {\bar f_1^{5}\over
f_1}\bar\partial_1^{10}(\phi),\\
& <X_2,\phi>=\int {\bar R_1(f_1,f_2)\over R_1(f_1,f_2)
}\bar\partial_2^{4}\big({\bar\partial(\bar f_1^{5})\over f_1}\bar
\partial^{10}(\phi)\big).\end{cases}$$

\bigskip

Nous allons donner un autre exemple d'une intersection compl\`ete au
voisinage de $0$ dans $\C^3$ :
$$(f_1,f_2,f_3):=(z_1^3+z_2^3+z_3^3,z_1^2+z_2^2,z_1^2+z_2z_3);$$
comme pr\'ec\'edement des calculs \'el\'ementaires conduisent aux
r\'esultants :
$$R(f_1,f_2)=2z_2^6+2z_2^3z_3^3+z_3^6,$$
$$R(f_1,f_3)=z_2^6+3z_2^3z_3^3+z_3^6;$$
De la m\^eme mani\`ere, on obtient que le r\'esultant de ces deux
derni\`eres expressions est $z_3^9$. Ainsi, il existe une matrice
$A$(que l'on peut calculer avec l'algoritme d'Euclide), qui
transforme $(f_1,f_2,f_3)$ en l'intersection compl\`ete
$(z_1^3+z_2^3+z_3^3,2z_2^6+2z_2^3z_3^3+z_3^6,
z_3^9):=(g_1,g_2,g_3)$. Pour cette derni\`ere, il est facile de
construire la famille de courants d\'ecrivant l'id\'eal en $0$
engendr\'e par $(g_1,g_2,g_3)$, en choisissant $N_1$ et $N_2$ deux
entiers suffisamments grands afin que les expressions ci-dessous
soient int\'egrables :
$$\begin{cases} & <X_1,\phi>=\int{\bar g_1^{N_1}\over
g_1}\bar\partial_1^{3N_1}(\phi),\\
& <X_2,\phi>=\int{\bar g_2^{N_2}\over
g_2}\bar\partial_2^{6N_2}\big({\bar\partial (\bar g_1^{N_1})\over
g_1}\bar\partial_1^{3N_1}(\phi)\big),\\
& <X_3,\phi>=\int{\bar g_3\over
g_3}\bar\partial_3^{9}\big({\bar\partial(\bar g_2)^{N_2}\over
g_2}\bar\partial_2^{6N_2}\big({\bar\partial(\bar g_1)^{N_1}\over
g_1}\bar\partial_1^{3N_1}(\phi)\big)\big).\end{cases}$$En utilisant
les notations de la section $3$, la famille de courants
$(X_1,\delta_2X_2,\delta_3X_3)$ d\'ecrit donc l'id\'eal en $0$
engendr\'e par $(f_1,f_2,f_3)$.

\bigskip

Consid\'erons dans $\C^2$, l'intersection compl\`ete :
$$f=(z_2^2+z_1^2+z_1^3,z_1^{4})=(f_1,f_2).$$ Cette exemple est
important car Tsikh-Passare, dans [12], montre que pour obtenir le
courant r\'esiduel associ\'e \`a $f$, il est n\'ecessaire de
consid\'erer $\theta\rightarrow\varepsilon(\theta)$ un pav\'e
admissible (les limites d\'ependent de la mani\`ere de "tendre vers
z\'ero"). N\'eanmoins, nous remarquons que $f$ est d\'ej\`a du type
r\'egulier et il est ais\'e de voir que les courants suivants
d\'ecrivent l'id\'eal local de $f$ :
$$\begin{cases} & <X_1,\phi>=\int{\bar f_1^5\over
f_1}\bar\partial_2^{10}(\phi),\\
& <X_2,\phi>=\int{\bar z_1^4\over
z_1^4}\bar\partial_1^4\big({\bar\partial(\bar f_1^5)\over
f_1}\bar\partial_2^{10}(\phi)\big).\end{cases}$$

\end{document}